\overfullrule=0pt
\centerline {\bf A more complete version of a minimax theorem}\par
\bigskip
\bigskip
\centerline {BIAGIO RICCERI}\par
\bigskip
\bigskip
{\bf Abstract.} In this paper, we present a more complete version of the minimax theorem established in [7]. As a consequence, we get,
for instance, the following result: Let $X$ be a compact, not singleton subset of a normed space $(E,\|\cdot\|)$ and let $Y$ be a convex subset
of $E$ such that $X\subseteq \overline {Y}$.
Then, for every convex set $S\subseteq Y$ dense in $Y$, for every upper semicontinuous bounded function $\gamma:X\to {\bf R}$ and
for every $\lambda>{{4\sup_X|\gamma|}\over {\hbox {\rm diam}(X)}}$,
there exists $y^*\in S$ such that the function $x\to \gamma(x)+\lambda\|x-y^*\|$ has at least two global maxima in $X$.\par
\par
\bigskip
\bigskip
Here and in what follows, $X$ is a topological space and $Y$ is a convex set in a real Hausdorff topological vector space.
\smallskip
A function $h:X\to {\bf R}$ is said to be inf-compact if $h^{-1}(]-\infty,r])$ is compact for all $r\in {\bf R}$. \par
\smallskip
A function $k:Y\to {\bf R}$ is said to be quasi-concave (resp. quasi-convex)) $k^{-1}([r,+\infty[)$
(resp. if $k^{-1}(]-\infty,r])$ is convex for all $r\in {\bf R}$.\par
\smallskip
If $S$ is a convex subset of $Y$, we denote by ${\cal A}_S$ the class of all functions $f:X\times Y\to {\bf R}$ such that, for each $y\in S$, the function $f(\cdot,y)$ is lower semicontinuous and $\inf$-compact. 
\par
\smallskip
 Moreover, we denote by ${\cal B}$ the class of all functions $f:X\times Y\to {\bf R}$ such that
either, for each $x\in X$, the function $f(x,\cdot)$ is quasi-concave and continuous, or, for each $x\in X$, the function $f(x,\cdot)$ is concave.\par
\smallskip
 For any $f:X\times Y\to {\bf R}$, we set
$$\alpha_f=\sup_Y\inf_Xf$$
and
$$\beta_f=\inf_X\sup_Yf\ .$$
Also, we denote by ${\cal C}_f$ the family of all sets $S\subseteq Y$ such that
$$\inf_X\sup_Sf=\inf_X\sup_Yf$$
and by $\tilde{\cal C}_f$ the family of all sets $S\subseteq Y$ such that
$$\sup_{y\in S}f(x,y)=\sup_{y\in Y}f(x,y)$$
for all $x\in X$.\par
\smallskip
In particular, notice that $S\in \tilde{\cal C}_f$ provided, for each $x\in X$, there is a topology on $Y$ for which $S$ is dense and $f(x,\cdot)$ is
lower semicontinuous.\par
\smallskip
Furthermore, we denote by $\tau_f$ the topology on $Y$ generated by the family
$$\{\{y\in Y : f(x,y)< r\}\}_{x\in X, r\in {\bf R}}\ .$$
So, $\tau_f$ is the weakest topology on $Y$ for which $f(x,\cdot)$ is upper semicontinuous 
for all $x\in X$. 
\smallskip
In [7], we established the following minimax result:\par
\medskip
THEOREM A. - {\it For every $g\in {\cal A}_Y\cap {\cal B}$, at least one of the following assertions holds:\par
\noindent
$(j)$\hskip 5pt $\sup_Y\inf_Xg=\inf_X\sup_Yg$\ ;\par
\noindent
$(jj)$\hskip 5pt there exists $y^*\in Y$ such that the function $g(\cdot,y^*)$ has at least two global minima.}\par
\medskip
The relevance of Theorem A resides essentially in the fact that it is a flexible tool which can fruitfully be used to obtain
meaningful results of various nature. This is clearly shown by a series of recent papers ([8]-[14]).\par
\smallskip
So, we believe that it is of interest to present a more complete form of Theorem A: this is just the aim of this paper.\par
\smallskip
Here is the main abstract result (with the usual rules in $\overline {\bf R}$):
\medskip
THEOREM 1. - {\it Let $f:X\times Y\to {\bf R}$. Assume that there is a function $\psi:Y\to {\bf R}$ such that $f+\psi\in {\cal B}$
and
$$\alpha_{f+\psi}<\beta_{f+\psi}\ .$$
Then, for every convex set $S\in {\cal C}_{f+\psi}$, for every bounded function $\varphi:X\to {\bf R}$  and 
for every $\lambda>0$ such that $\lambda f+\varphi\in {\cal A}_S$ and
$$\lambda>{{2\sup_X|\varphi|}\over {\beta_{f+\psi}-\alpha_{f+\psi}}}\ ,\eqno{(1)}$$
there exists $y^*\in S$ such that the function $\lambda f(\cdot,y^*)+\varphi(\cdot)$ has at least two global minima.}\par
\smallskip
PROOF. Consider the function $g:X\times Y\to {\bf R}$ defined by
$$g(x,y)=\lambda(f(x,y)+\psi(y))+\varphi(x)$$
for all $(x,y)\in X\times Y$. Since $S\in {\cal C}_{f+\psi}$, we have
$$\inf_X\sup_S(f+\psi)=\inf_X\sup_Y(f+\psi)\ . \eqno{(2)}$$
 So, taking $(1)$ and $(2)$ into account, we have
$$\sup_S\inf_Xg\leq \sup_Y\inf_Xg\leq \lambda\alpha_{f+\psi}+\sup_{X}|\varphi|$$
$$<\lambda\beta_{f+\psi}-\sup_X|\varphi|=\lambda\inf_X\sup_S(f+\psi)-\sup_X|\varphi|
\leq \inf_X\sup_Sg\ .\eqno{(3)}$$
Now, observe that $g\in {\cal A}_S$ since $\lambda f+\varphi\in {\cal A}_S$ and, at the same time, $g\in {\cal B}$ since
$f+\psi\in {\cal B}$. As a consequence, we can apply Theorem A to the restriction of the function $g$
to $X\times S$. Therefore, in view of $(3)$, there exists of $y^*\in S$ such that 
the function $g(\cdot,y^*)$ (and hence $\lambda f(\cdot,y^*)+\varphi(\cdot))$ has at least two global minima, as claimed. 
\hfill $\bigtriangleup$\par
\medskip
REMARK 1. As the above proof shows, Theorem 1 is a direct consequence of Theorem A. However, there are essentially four advantages
of Theorem 1 with respect to Theorem A. Namely, suppose that, for a given function $g\in {\cal A}_Y$, we are interested in
ensuring the validity of assertion $(jj)$. Then, if we apply Theorem A in this regard, we have to show that $g\in {\cal B}$ and that assertion $(j)$ does not hold. To the contrary, if we apply Theorem 1, we can ensure the validity of $(jj)$ also in cases where either $g\not\in {\cal B}$ or
$(j)$ holds true too. In addition, Theorem 1 is able to ensure the validity of $(jj)$ even in a remarkably stronger way: not only extending it to suitable perturbations of $g$, but also offering an information on the location of $y^*$.\par
\medskip
First, we wish to show how to obtain the very classical minimax theorems in [3] and [6] by means of Theorem 1.\par
\smallskip
Let $V$ be a real vector space, $A\subseteq V$, $\varphi:A\to {\bf R}$. We say that $\varphi$ is finitely lower semicontinuous if, for every
finite-dimensional linear subspace $F\subseteq V$, the function $f_{|A\cap F}$ is lower semicontinuous in the Euclidean topology of $F$.\par
\medskip
In the next result, the topology of $X$ has no role.\par
\medskip
THEOREM 2. - {\it Let $X$ be a convex set in a real vector space and let $f\in {\cal B}$. Assume that there is a convex set
$S\in \tilde{\cal C}_f$ such that $f(\cdot,y)$ is finitely lower semicontinous and convex for all $y\in S$. Finally,
assume that, for some $x_0\in X$, the function that $f(x_0,\cdot)$ is $\tau_f-\sup$-compact.\par
Then, one has
$$\sup_Y\inf_Xf=\inf_X\sup_Yf\ .$$}\par
\smallskip
PROOF. Arguing by contradiction, assume that
$$\sup_Y\inf_Xf<\inf_X\sup_Yf\ .$$
Denote by ${\cal D}$ the family of all convex polytopes in $X$. Since ${\cal D}$ is a filtering cover of $X$ and $f(x_0,\cdot)$ is $\tau_f-\sup$-compact,
by Proposition 2.1 of [7], there exists $P\in {\cal D}$ such that
$$\sup_Y\inf_Pf<\inf_P\sup_Yf\ .$$ Let $\|\cdot\|$ be the Euclidean norm on
span$(P)$. So, $\|\cdot\|^2$ is strictly convex. Now, fix $\lambda$ so that
$$\lambda>{{2\sup_{x\in P}\|x\|^2}\over
{\inf_P\sup_Yf-\sup_Y\inf_Pf}}\ .$$ Notice that, for each $y\in S$, the function
$x\to \|x\|^2+\lambda f(x,y)$ is inf-compact in $P$ with respect to the Euclidean topology. As a consequence, if we consider
$P$ equipped with the Euclidean topology, we can apply Theorem 1 to the restriction of $f$ to $P\times Y$ (recall that $S\in \tilde{\cal C}_f$), taking
$\varphi=\|\cdot\|^2$. Accordingly, there
would exist $y^*\in S$ such that the function $x\to \|x\|^2+\lambda f(x,y)$ has at least two global minima in $P$. But, this is absurd since
this function is strictly convex.\hfill $\bigtriangleup$
\medskip
Reasoning exactly as in the proof of Theorem 2 (even in a simplified way, since there is no need to consider the family ${\cal D}$), we also get\par
\medskip
THEOREM 3. - {\it Let $X$ be a compact convex set in a topological vector space such that there exists a lower semicontinuous, strictly
convex, bounded function $\varphi:X\to {\bf R}$. Let $f\in {\cal B}$.  Assume that there is a convex set
$S\in {\cal C}_f$ such that $f(\cdot,y)$ is lower semicontinuous and convex for all $y\in S$.\par
Then, one has
$$\sup_Y\inf_Xf=\inf_X\sup_Yf\ .$$}\par
\medskip
We now revisit two applications of Theorem A in the light of Theorem 1.\par
\smallskip
The first one concerns the so called farthest points ([1], [4]).\par
\medskip
THEOREM 4. - {\it Let $X$ be compact, not singleton and let $(E,d)$ be a metric space such that $X\subseteq E$. Moreover, let $h:Y\to E$ be such that $X\subseteq \overline {h(Y)}$ and let the function $(x,y)\to f(x,y):=-d(x,h(y))$ belong to ${\cal B}$.\par
Then, for every convex set $S\in {\cal C}_f$, for every bounded function $\gamma:X\to {\bf R}$ and for every $\lambda>0$
 such that $\lambda f-\gamma\in {\cal A}_S$  and
$$\lambda>{{4\sup_X|\gamma|}\over {\hbox {\rm diam}(X)}}\ ,$$
there exists $y^*\in S$ such that the function $x\to \gamma(x)+\lambda d(x,h(y^*))$ has at least two global maxima in $X$.}\par
\smallskip
PROOF. Since $X\subseteq \overline {h(Y)}$, we have
$$\sup_{x\in X}\inf_{y\in Y}d(x,h(y))=0\ .\eqno{(4)}$$
Also, for each $x_1, x_2\in X$, $y\in Y$, we have
$${{d(x_1,x_2)}\over {2}}\leq \max\{d(x_1,h(y)), d(x_2,h(y))\}$$
and so
$${{\hbox {\rm diam}(X)}\over {2}}\leq\inf_{y\in Y}\sup_{x\in X}d(x,h(y))\ .\eqno{(5)}$$
Hence, in view of $(4)$ and $(5)$, we have
$$\sup_Y\inf_Xf\leq -{{\hbox {\rm diam}(X)}\over {2}}<0=\inf_X\sup_Yf\ .$$
Now, the conclusion follows directly from Theorem 1 taking $\varphi=-\gamma$.\hfill $\bigtriangleup$\par
\medskip
Of course, the most natural corollary of Theorem 4 is as follows:\par
\medskip
COROLLARY 1. - {\it Let $X$ be a compact, not singleton subset of a normed space $(E,\|\cdot\|)$ and let $Y$ be a convex subset
of $E$ such that $X\subseteq \overline {Y}$.\par
Then, for every convex set $S\subseteq Y$ dense in $Y$, for every upper semicontinuous bounded function $\gamma:X\to {\bf R}$ and
for every $\lambda>{{4\sup_X|\gamma|}\over {\hbox {\rm diam}(X)}}\ ,$
there exists $y^*\in S$ such that the function $x\to \gamma(x)+\lambda\|x-y^*\|$ has at least two global maxima in $X$.}\par
\medskip
In turn, from Corollary 1, we clearly get\par
\medskip
COROLLARY 2. - {\it Let $X$ be a compact subset of a normed space $(E,\|\cdot\|)$ and let $Y$ be a convex subset
of $E$ such that $X\subseteq \overline {Y}$. Assume that there exist a sequence $\{S_n\}$ of convex subsets of $Y$ dense in $Y$ and
a sequence $\{\gamma_n\}$ of upper semicontinuous bounded real-valued functions on $X$, with $\lim_{n\to \infty}\sup_X|\gamma_n|=0$,
such that, for each $n\in {\bf N}$ and for each $y\in S_n$, the function $x\to \gamma_n(x)+\|x-y\|$ has a unique global maximum in $X$.\par
Then, $X$ is a singleton.}\par
\medskip
REMARK 2. - Notice that Corollary 2 improves Theorem 1.1 of [14] which, in turn, extended a classical result by Klee ([5]) to normed spaces. More
precisely, Theorem 1.1 of [14] agrees with the particular case of Corollary 2 in which each $S_n$ is equal to conv$(X)$ and each $\gamma_n$ is
equal to $0$.
\medskip
The second application concerns the calculus of variations.\par
\smallskip
Let $\Omega\subset {\bf R}^n$ be a bounded domain with smooth boundary and let $p>1$. On the Sobolev space $W^{1,p}(\Omega)$, we
consider the norm
$$\|u\|=\left ( \int_{\Omega}|\nabla u(x)|^p dx+\int_{\Omega}|u(x)|^p dx\right ) ^{1\over p}\ .$$
If $n\geq p$, we denote by ${\cal E}$ the class of all
continuous functions $\sigma:{\bf R}\to {\bf R}$ such that
$$\sup_{\xi\in {\bf R}}{{|\sigma(\xi)|}\over
{1+|\xi|^q}}<+\infty\ ,$$
where  $0<q< {{pn}\over {n-p}}$ if $p<n$ and $0<q<+\infty$ if
$p=n$. While, when $n<p$, ${\cal E}$ stands for the class
of all continuous functions $\sigma:{\bf R}\to {\bf R}$.\par
\smallskip
Recall that a function $h:\Omega\times {\bf R}^{m}\to {\bf R}$
 is said to be a normal integrand ([15]) if it is
${\cal L}(\Omega)\otimes {\cal B}({\bf R}^m)$-measurable and
$h(x,\cdot)$ is lower semicontinuous for a.e. $x\in \Omega$.
Here ${\cal L}(\Omega)$ and ${\cal B}({\bf R}^m)$ denote the
Lebesgue and the Borel $\sigma$-algebras of subsets of $\Omega$ and
${\bf R}^m$, respectively. \par
\smallskip
Recall that if $h$ is a normal integrand then, for each measurable function $u:\Omega\to
{\bf R}^m$, the composite function $x\to h(x,u(x))$ is measurable
([15]).\par
\smallskip
We denote by ${\cal F}$ the class of all normal integrands  $h:\Omega\times {\bf R}\times {\bf R}^{n}\to {\bf R}$
such that $h(x,\xi,\cdot)$ is convex for all $(x,\xi)\in \Omega\times {\bf R}$ and there are $M\in L^1(\Omega)$, $b>0$ such that
$$M(x)-b(|\xi|+|\eta|)|)\leq h(x,\xi,\eta)$$
for all $(x,\xi,\eta)\in \Omega\times {\bf R}\times {\bf R}^{n}\to {\bf R}$.\par
\smallskip
Let us also recall two results proved in [9].\par
\medskip
PROPOSITION 1. - {\it Let $\Omega\subset {\bf R}^n$ be a bounded domain with smooth boundary, let $p>1$ and let
$h:\Omega\times {\bf R}\times {\bf R}^n\to {\bf R}$ be normal integrand such that, for some $c, d>0$, one has
$$c|\eta|^p-d\leq h(x,\xi,\eta)$$
for all $(x,\xi,\eta)\in\Omega\times{\bf R}\times{\bf R}^n$ and
$$\lim_{|\xi|\to +\infty}\inf_{(x,\eta)\in \Omega\times{\bf R}^n}h(x,\xi,\eta)=+\infty\ .$$
Then, in $W^{1,p}(\Omega)$, one has
$$\lim_{\|u\|\to +\infty}\int_{\Omega}h(x,u(x),\nabla u(x))dx=+\infty\ .$$}\par
\medskip
PROPOSITION 2. - {\it Let $X, Y$ be two non-empty sets and $I:X\to {\bf R}$, $J:X\times Y\to {\bf R}$ two given
functions. Assume that there are two sets $A, B\subset X$ such that:\par
\noindent
$(a)$\hskip 5pt $\sup_AI<\inf_BI$\ ;\par
\noindent
$(b)$\hskip 5pt $\sup_Y\inf_AJ(x,y)\leq 0$\ ;\par
\noindent
$(c)$\hskip 5pt $\inf_B\sup_YJ(x,y)\geq 0$\ ;\par
\noindent
$(d)$\hskip 5pt $\inf_{X\setminus B}\sup_YJ(x,y)=+\infty$\ .\par
Then, one has
$$\sup_Y\inf_X(I+J)\leq \sup_AI<\inf_BI\leq\inf_X\sup_Y(I+J)\ .$$}\par
\medskip
Furthermore, let us also recall the following classical fact:\par
\medskip
PROPOSITION 3. - {\it Let $A\subseteq {\bf R}^n$ be any open set and let $v\in L^1(A)\setminus
\{0\}$.\par
Then, one has
$$\sup_{\alpha\in C^{\infty}_0(A)}\int_{A}\alpha(x)v(x)dx=+\infty\ .$$}
\medskip
After these preliminaries, we can prove the following result:\par
\medskip
THEOREM 5. - {\it  Let $h, k\in {\cal F}$ and
let $\sigma\in {\cal E}$ be a strictly monotone function. Assume that:\par
\noindent
$(i)$\hskip 5pt there are $c, d>0$ such that
$$c|\eta|^p-d\leq h(x,\xi,\eta)$$
for all $(x,\xi,\eta)\in \Omega\times{\bf R}\times{\bf R}^n$ and
$$\lim_{|\xi|\to +\infty}{{\inf_{(x,\eta)\in \Omega\times{\bf R}^n}h(x,\xi,\eta)}\over {|\sigma(\xi)|+1}}=+\infty\ ;$$
\noindent
$(ii)$\hskip 5pt for each $\xi\in {\bf R}$, the function $h(\cdot,\xi,0)$ lies in $L^1(\Omega)$\ ;\par
\noindent
$(iii)$\hskip 5pt there are $\xi_1, \xi_2, \xi_3\in {\bf R}$, with $\xi_1<\xi_2<\xi_3$, such that
$$\max\left \{\int_{\Omega}h(x,\xi_1,0)dx, \int_{\Omega}h(x,\xi_3,0)dx\right \}<\int_{\Omega}h(x,\xi_2,0)dx\ .$$
Then, for every sequentially weakly closed set $V\subseteq W^{1,p}(\Omega)$, containing the constants, for every convex set
$T\subseteq L^{\infty}(\Omega)$ dense in $L^{\infty}(\Omega)$, for every non-decreasing, continuous, bounded function
$\omega:U\to {\bf R}$, where $U:=\{\int_{\Omega}k(x,u(x),\nabla u(x))dx : u\in W^{1,p}(\Omega)\}$,
and for every $\lambda$ satisfying
$$\lambda>
{{2\sup_U|\omega|}\over {\int_{\Omega}h(x,\xi_2,0)dx-\max\left \{\int_{\Omega}h(x,\xi_1,0)dx, \int_{\Omega}h(x,\xi_3,0)dx\right \}}}\ ,\eqno{(6)}$$
there exists $\gamma\in T$
such that the restriction to $V$ of the functional 
$$u\to \lambda\int_{\Omega}h(x,u(x),\nabla u(x))dx+\int_{\Omega}\gamma(x)\sigma(u(x))dx
+\omega\left ( \int_{\Omega}k(x,u(x),\nabla u(x))dx\right )$$
has at least two global minima. The same conclusion holds also with $T=C^{\infty}_0(\Omega)$.}\par
\smallskip
PROOF. Fix $V, T, \omega, \lambda$ as in the conclusion.
Since $\sigma\in {\cal E}$, in view of the Rellich-Kondrachov theorem, for each $u\in W^{1,p}(\Omega)$, we have $\sigma\circ u\in L^1(\Omega)$ and, for each $\gamma\in L^{\infty}(\Omega)$,
the functional $u\to \int_{\Omega}\gamma(x)\sigma(u(x))dx$ is sequentially weakly continuous. Moreover, since $h, k\in {\cal F}$
the functionals $u\to \int_{\Omega}h(x,u(x),\nabla u(x)dx$ and $u\to \int_{\Omega}k(x,u(x),\nabla u(x)dx$ (possibly taking the value $+\infty$)
are sequentially weakly lower semicontinuous ([2], Theorem 4.6.8). Hence, since $\omega$ is non-decreasing and continuous, the functional
$u\to \omega\left (\int_{\Omega}k(x,u(x),\nabla u(x)dx\right )$ is sequentailly weakly lower semicontinuous too.
Set
$$X=\left \{u\in V : \int_{\Omega}h(x,u(x),\nabla u(x))dx<+\infty\right \}\ .$$
By $(ii)$, the constants belong to $X$.  
Fix $\gamma\in L^{\infty}(\Omega)$. By $(i)$, there is $\delta>0$ such that
$$h(x,\xi,\eta)-2\|\gamma\|_{L^{\infty}(\Omega)}|\sigma(\xi)|\geq 0$$
for all $(x,\xi,\eta)\in \Omega\times{\bf R}\times{\bf R}^n$ with $|\xi|>\delta$. So, we have
$${{c}\over {2}}|\eta|^p-d-\|\gamma\|_{L^{\infty}(\Omega)}\sup_{|\xi|\leq\delta}|\sigma(\xi)|\leq h(x,\xi,\eta)+\gamma(x)\sigma(\xi)$$
for all $(x,\xi,\eta)\in \Omega\times{\bf R}\times{\bf R}^n$ and, of course,
$$\lim_{|\xi|\to +\infty}\inf_{(x,\eta)\in \Omega\times{\bf R}^n}(h(x,\xi,\eta)+\gamma(x)\sigma(\xi))=+\infty\ .$$
Consequently, in view of Proposition 2.1, we have, in $W^{1,p}(\Omega)$,
$$\lim_{\|u\|\to +\infty}\left (\int_{\Omega}h(x,u(x),\nabla u(x))dx+\int_{\Omega}\gamma(x)\sigma(u(x))dx\right )=+\infty \ .$$
This implies that, for each $r\in {\bf R}$, the set
$$\left \{u\in V: \int_{\Omega}h(x,u(x),\nabla u(x))dx+\int_{\Omega}\gamma(x)\sigma(u(x))dx\leq r\right \}$$
is weakly compact by reflexivity and by Eberlein-Smulyan's theorem. Of course, we also have
$$\left \{u\in V: \int_{\Omega}h(x,u(x),\nabla u(x))dx+\int_{\Omega}\gamma(x)\sigma(u(x))dx\leq r\right \}$$
$$=\left \{u\in X: \int_{\Omega}h(x,u(x),\nabla u(x))dx+\int_{\Omega}\gamma(x)\sigma(u(x))dx\leq r\right \}\ .$$
Now, observe that, if we put
$$A=\{\xi_1, \xi_3\}$$
and
$$B=\{\xi_2\}\ ,$$
and define $I:X\to {\bf R}$, $J:X\times L^{\infty}(\Omega)\to {\bf R}$ by
$$I(u)=\int_{\Omega}h(x,u(x),\nabla u(x))dx \ ,$$
$$J(u,\gamma)=\int_{\Omega}\gamma(x)(\sigma(u(x))-\sigma(\xi_2))dx$$
for all $u\in X$, $\gamma\in L^{\infty}(\Omega)$, we clearly have
$$\inf_{u\in B}\sup_{\gamma\in L^{\infty}(\Omega)}J(u,\gamma)=0$$
and, by $(iii)$,
$$\sup_AI<\inf_BI\ .$$
Since $\sigma$ is strictly monotone, the numbers $\sigma(\xi_1)-\sigma(\xi_2)$ and $\sigma(\xi_3)-\sigma(\xi_2)$
have opposite signs. This clearly implies that
$$\sup_{\gamma\in L^{\infty}(\Omega)}\inf_{u\in A}J(u,\gamma)\leq 0\ .$$
Furthermore, if $u\in X\setminus \{\xi_2\}$, again by strict monotonicity, $\sigma\circ u\neq \sigma(\xi_2)$, and so we have
$$\sup_{\gamma\in L^{\infty}(\Omega)}J(u,\gamma)=+\infty\ .$$
Therefore, the sets $A, B$ and the functions $I, J$ satisfy the assumptions of Proposition 2 and hence we have
$$\sup_{L^{\infty}(\Omega)}\inf_X(I+J)\leq\max\left \{\int_{\Omega}h(x,\xi_1,0)dx, \int_{\Omega}h(x,\xi_3,0)dx\right \}<\int_{\Omega}h(x,\xi_2,0)dx\leq
\inf_X\sup_{L^{\infty}(\Omega)}(I+J)\ .\eqno{(7)}$$
Now, we can apply Theorem 1 considering $X$ equipped with the weak topology and taking
$$Y=L^{\infty}(\Omega)\ ,$$
$$f=I+J\ ,$$
$$\psi=0\ ,$$
$$S={{1}\over {\lambda}}T$$
and
$$\varphi(u)=\omega\left ( \int_{\Omega}k(x,u(x),\nabla u(x))dx\right )\ .$$
Notice that, in view of $(7)$, inequality $(1)$ holds thanks to $(6)$, and the conclusion follows.
When $T=C^{\infty}_0(\Omega)$ the same proof as above holds in view of Proposition 3.
\hfill $\bigtriangleup$\par
\medskip
REMARK 4. - Notice that condition $(iii)$ holds if and only if the function $\int_{\Omega}h(x,\cdot,0)$ is not quasi-convex.\par
\medskip
REMARK 5. - For $\omega=0$, Theorem 5 reduces to Theorem 1.2 of [9].\par
\medskip
 We conclude presenting an application of Theorem 5 to the Neumann problem for a Kirchhoff-type equation.\par
\smallskip
Given $K:[0,+\infty[\to {\bf R}$ and a Carath\'eodory function $\psi:\Omega\times {\bf R}\to {\bf R}$, consider the following Neumann problem
$$\cases {-K\left( \int_{\Omega}|\nabla u(x)|^pdx\right )\hbox {\rm div}(|\nabla u|^{p-2}\nabla u)=
\psi(x,u)
 & in
$\Omega$\cr & \cr {{\partial u}\over {\partial \nu}}=0 & on
$\partial \Omega$\ ,\cr} $$
where $\nu$ is the outward unit normal to $\partial\Omega$.\par
 \par
Let us recall
that a weak solution
of this problem is any $u\in W^{1,p}(\Omega)$ such that, for every $v\in W^{1,p}(\Omega)$, one has $\psi(\cdot,u(\cdot))v(\cdot)\in L^1(\Omega)$
and
 $$K\left( \int_{\Omega}|\nabla u(x)|^pdx\right )\int_{\Omega}|\nabla u(x)|^{p-2}\nabla u(x)\nabla v(x)dx
-\int_{\Omega}
\psi(x,u(x))v(x)dx=0\ .$$
\medskip
THEOREM 6. - {\it Let $f, g:{\bf R}\to {\bf R}$ be two $C^1$ functions lying in ${\cal E}$ and satisfying the following
conditions:\par
\noindent
$(a_1)$\hskip 5pt the function $g'$ has a constant sign and $\hbox {\rm int}((g')^{-1}(0))=\emptyset$\ ;\par
\noindent
$(a_2)$\hskip 5pt $\lim_{|\xi|\to +\infty}{{f(\xi)}\over {|g(\xi)|+1}}=+\infty$\ ;\par
\noindent
$(a_3)$\hskip 5pt there are $\xi_1, \xi_2, \xi_3\in {\bf R}$, with $\xi_1<\xi_2<\xi_3$, such that
$$\max\{f(\xi_1), f(\xi_3)\}<f(\xi_2)\ .$$
Then, for every $a>0$, for every $\beta\in L^{\infty}(\Omega)$, with $\inf_{\Omega}\beta>0$,  for every convex set $T\subseteq L^{\infty}(\Omega)$ dense
in $L^{\infty}(\Omega)$, for every $C^1$, non-decreasing, bounded function $\chi:[0,+\infty[\to {\bf R}$, and for every $\lambda$
satisfying
$$\lambda>{{2\sup_{[0,+\infty[}|\chi|}\over {p(f(\xi_2)-\max\{f(\xi_1), f(\xi_3)\})\int_{\Omega}\beta(x)dx}}$$
there exists $\gamma\in T$ such that the problem
$$\cases {-\left ( a+\chi'\left ( \int_{\Omega}|\nabla u(x)|^pdx\right )\right )
\hbox {\rm div}(|\nabla u|^{p-2}\nabla u)=\gamma(x)g'(u)-\lambda\beta(x)f'(u)
 & in
$\Omega$\cr & \cr {{\partial u}\over {\partial \nu}}=0 & on
$\partial \Omega$\ ,\cr}\eqno{(P)}$$
has at least two weak solutions.}\par
\smallskip
PROOF. Fix $a, \beta$, $T$, $\chi$ and $\lambda$ as in the conclusion. We are going to apply Theorem 5, defining $h, k, \sigma$ by
$$h(x,\xi,\eta)={{a}\over {p\lambda}}|\eta|^p+\beta(x)f(\xi)\ ,$$
$$k(\eta)=|\eta|^p\ ,$$
$$\sigma(\xi)=-g(\xi)$$
for all $(x,\xi,\eta)\in \Omega\times {\bf R}\times {\bf R}^n$. It is immediate to realize that, by $(a_1)-(a_3)$, the above $h, k, \sigma$
satisfy the assumptions of Theorem 5. Then, applying Theorem 5 with
$\omega={{1}\over {p}}\chi$, we get the existence of $\gamma\in T$ such that the functional
$$u\to \lambda\left ({{a}\over {p\lambda}}\int_{\Omega}|\nabla u(x)|^pdx+\int_{\Omega}\beta(x)f(u(x))dx\right )
-\int_{\Omega}\gamma(x)g(u(x))dx+{{1}\over {p}}\chi\left ( \int_{\Omega}\nabla u(x)|^pdx\right )$$
has at least two global minima in $W^{1,p}(\Omega)$. But, by classical results (recall that $f, g\in {\cal E}$),
 such a functional is $C^1$ and its critical points (and
so, in particular, its global minima) are 
weak solutions of problem $(P)$. The proof is complete.\hfill $\bigtriangleup$
\medskip
A challenging problem is as follows:\par
\medskip
PROBLEM 1. - Does the conclusion of Theorem 6 hold with {\it three} instead of {\it two} ?
\vfill\eject
\centerline {\bf References}\par
\bigskip
\bigskip
\noindent
[1]\hskip 5pt S. COBZA\c{S}, {\it Geometric properties of Banach spaces and the existence of nearest and farthest points},
Abstr. Appl. Anal., 2005, n. {\bf 3}, 259-285.\par
\smallskip
\noindent
[2]\hskip 5pt Z. DENKOWSKI, S. MIG\'ORSKI and N. S. PAPAGEORGIOU,
{\it An Introduction to Nonlinear Analysis: Applications}, Kluwer Academic
Publishers, 2003.\par
\smallskip
\noindent
[3]\hskip 5pt K. FAN, {\it Minimax theorems}, Proc. Nat. Acad. Sci. U.S.A., {\bf 39} (1953), 42-47.\par
\smallskip
\noindent
[4]\hskip 5pt J.-B. HIRIART-URRUTY, {\it La conjecture des points les plus \'eloign\'es revisit\'ee}, Ann. Sci. Math. Qu\'ebec {\bf 29} (2005), 197-214.
\par
\smallskip
\noindent
[5]\hskip 5pt V. L. KLEE, {\it Convexity of Chebyshev sets}, Math. Ann.,
{\bf 142} (1960/1961), 292-304.\par
\smallskip
\noindent
[6]\hskip 5pt H. KNESER, {\it Sur un th\'eor\`eme fondamental de la th\'eorie des jeux,} C. R. Acad. Sci. Paris {\bf 234} (1952), 2418-2420.\par
\smallskip
\noindent
[7]\hskip 5pt B. RICCERI, {\it On a minimax theorem: an improvement, a new proof and an overview of its applications},
Minimax Theory Appl., {\bf 2} (2017), 99-152.\par
\smallskip
\noindent
[8]\hskip 5pt B. RICCERI, {\it Another multiplicity result for the periodic solutions of certain systems}, Linear Nonlinear Anal., {\bf 5} (2019), 371-378.\par
\smallskip
\noindent
[9]\hskip 5pt B. RICCERI, {\it Miscellaneous applications of certain minimax theorems II}, Acta Math. Vietnam., {\bf 45} (2020), 515-524.\par
\smallskip
\noindent
[10]\hskip 5pt B. RICCERI, {\it A class of equations with three solutions}, Mathematics (2020), {\bf 8}, 478.\par
\smallskip
\noindent
[11]\hskip 5pt B. RICCERI, {\it An invitation to the study of a uniqueness problem}, in ``Nonlinear Analysis and Global Optimization", Th. M. Rassias and P. M. Pardalos eds., 445-448, Springer, 2021.\par
\smallskip
\noindent
[12]\hskip 5pt B. RICCERI, {\it A class of functionals possessing multiple global minima}, Stud. Univ. Babe\c{s}-Bolyai Math., {\bf 66} (2021), 75-84.\par
\smallskip
\noindent
[13]\hskip 5pt B. RICCERI, {\it An alternative theorem for gradient systems}, Pure Appl. Funct. Anal., {\bf 6} (2021), 373-381.\par
\smallskip
\noindent
[14]\hskip 5pt B. RICCERI, {\it On the applications of a minimax theorem}, Optimization, to appear.\par
\smallskip
\noindent
[15]\hskip 5pt R. T. ROCKAFELLAR, {\it  Integral functionals, normal integrands and measurable selections}, Lecture Notes in Math., Vol. 543,  157-207, Springer, Berlin, 1976.\par

\bigskip
\bigskip
\bigskip
\bigskip
Department of Mathematics and Informatics\par
University of Catania\par
Viale A. Doria 6\par
95125 Catania, Italy\par
{\it e-mail address}: ricceri@dmi.unict.it
\bye